\documentclass[11pt]{article}
\newtheorem{theorem}{Theorem}[section]
\newtheorem{corollary}{Corollary}[section]
\newtheorem{note}{Note}[section]
\usepackage{bm}
\usepackage{amssymb}
\usepackage{amsmath}
\usepackage{epsfig}

\begin{document}
\begin{center}
{\LARGE\bf The distribution of the maximum}\\[1ex]
{\LARGE\bf of an ARMA(1, 1) process}\\[1ex]
by\\[1ex]
Christopher S. Withers\\
Applied Mathematics Group\\
Industrial Research Limited\\
Lower Hutt, NEW ZEALAND\\[2ex]
Saralees Nadarajah\\
School of Mathematics\\
University of Manchester\\
Manchester M13 9PL, UK
\end{center}
\vspace{1.5cm}
{\bf Abstract:}~~We give the cumulative distribution function of $M_n$,  the maximum of a
sequence of $n$ observations from an ARMA(1, 1) process.
Solutions are first given in terms of repeated integrals
and then for the case, where the underlying random variables
are absolutely continuous.
The distribution of  $M_n$ is then given as a weighted sum of
the $n$th powers of the eigenvalues of a non-symmetric Fredholm kernel.
The weights are given in terms of the  left and right eigenfunctions of the kernel.

These results are large deviations expansions for estimates, since
the maximum need not be standardized to have a limit.
In fact, such a limit need not exist.


\section{Introduction and summary}
\setcounter{equation}{0}

There has been little work with respect to extreme value theory for ARMA (autoregressive-moving average) processes.
The authors are aware only of the work of Rootz\'{e}n (1983, 1986).
Both these papers provide the limiting extreme value distributions or assume that the errors come from a specific class.
We are aware of no work giving the {\it exact} distribution of the maximum of ARMA processes.

This paper continues the application of a powerful new method for
obtaining the {\it exact}
distribution of extremes of $n$ correlated observations as weighted
sums of $n$th powers of certain eigenvalues.
The method was first illustrated for a moving average of order 1 in Withers and Nadarajah (2009a)
and an autoregressive process of order 1  in Withers and Nadarajah (2009b).

Let  $\{e_i\}$ be independent and identically distributed random variables from some cumulative distribution function (cdf) $F$ on $R$.
Let $f (x)$ denote the probability density function (pdf) with respect to Lebesque measure.
We consider the ARMA process of order (1,1),
\begin{eqnarray*}
X_i - rX_{i-1}= e_i+ s e_{i-1}.
\end{eqnarray*}
We restrict ourselves to the most important case $r>0$.
(When this condition does not hold the method can be adapted as done in Withers and Nadarajah (2009b).)
In Section 2, we give expressions for the cdf of the maximum
\begin{eqnarray*}
M_n=\max_{i=1}^n X_i,
\
n\geq 1,
\end{eqnarray*}
in terms of repeated integrals.
This is obtained via the recurrence relationship
\begin{eqnarray}
G_{n}({\bf y}) = {\cal K}  G_{n-1}({\bf y}),
\
{\bf y} = (y_0,y_1),
\
n\geq 1,
\label{joint}
\end{eqnarray}
where
\begin{eqnarray}
G_{n}({\bf y}) = P(M_n\leq x,\ X_n\leq y_0,\ e_{n}\leq y_1),
\label{defG}
\end{eqnarray}
\begin{eqnarray}
{\cal K} h({\bf y})
&=&
E\ \int h( g_{\bf y}(s{\bf z} + e_0), d{\bf z}) I (e_{1}\leq y_1)
\nonumber
\\
&=&
r\int \int_{a_{\bf yz}}^\infty dz_0 f(y_{0x}-rz_0-sz_1)h(z_0,dz_1),
\label{calK}
\end{eqnarray}
\begin{eqnarray*}
g_{\bf y} ( t) = (y_{0x}-t)/r,
\
y_{0x}=\min(y_0,x),
\
a_{\bf yz}=( y_{0x}-sz_1-y_1)/r,
\end{eqnarray*}
$I(A)=1$ or 0 for $A$ true or false
and dependency on $x$ is suppressed.
So, $ {\cal K}$ is a linear  integral operator depending on $x$.
For (\ref{joint}) to work at $n=1$ we define $M_0=-\infty$ so  that
\begin{eqnarray}
G_{0}({\bf y}) = P(X_0\leq y_0,\ e_0\leq y_1).
\label{0}
\end{eqnarray}
Similarly,
\begin{eqnarray*}
G_{1}({\bf y}) = P(X_1\leq y_{0x},\ e_1\leq y_1)=G_0(y_{0x}, y_1).
\end{eqnarray*}
In Section 3, we consider the case when $F$ is absolutely continuous.
In this case we show that corresponding to ${\cal K}$ is a Fredholm kernel $K({\bf y}, {\bf z})$.
We give a solution in terms of its eigenvalues and eigenfunctions.
This leads easily to the asymptotic results stated in the abstract.
However, there are two problems: the kernel is a generalized function and
numerical solution
by direct Gaussian quadrature fails.
In Section 4, we show that these problems are avoided by using the iterated Fredholm kernel $K_2({\bf y}, {\bf z})$.

Our expansions for $P(M_n\leq x)$ for fixed $x$ are large deviation results.
If $x$ is replaced by $x_n$ such that  $P(M_n\leq x_n)$ tends to the generalized extreme value cdf,
then the expansion still holds, but not the asymptotic expansion
in terms of a single eigenvalue, since this may approach 1 as $n\rightarrow \infty$.

For $a$, $b$ functions on $R^2$, set $\int a=\int a({\bf y}) d{\bf y} = \int_{R^2} a({\bf y}) d{\bf y}$ and similarly for $\int ab$.

\section{Solutions using repeated integrals}
\setcounter{equation}{0}

\begin{theorem}
We have $G_n$ of (\ref{defG}) satisfying the recurrence relation (\ref{joint}) in terms
of the integral operator ${\cal K}$ of (\ref{calK}).
\end{theorem}

\noindent
{\bf Proof:}
For $n\geq 1$, $G_n$ of (\ref{defG}) satisfies
\begin{eqnarray*}
G_{n}({\bf y})
&=&
P(M_{n-1}\leq x,\ X_n\leq y_0,  rX_{n-1}+e_n +se_{n-1}\leq y_{0x},\ e_{n}\leq y_1)
\\
&=&
E\ P(M_{n-1}\leq x,  X_{n-1}\leq g_{\bf y} (se_{n-1}+e_n)|e_n)\ I( e_n\leq y_1)
\\
&=&
{\cal K}G_{n-1}({\bf y}).
\end{eqnarray*}
This ends the proof.
$\Box$\\

Our goal is to determine $u_n=P(M_n\leq x)= G_{n}(\pmb{\infty})$, where $\pmb{\infty}=(\infty, \infty)$.

\begin{theorem}
Set
\begin{eqnarray*}
a_n= [{\cal K}^n  G_0({\bf y})]_{{\bf y} = \pmb{\infty}},
\
n\geq 0,
\end{eqnarray*}
where $G_0$ is given by (\ref{0}).
Then
\begin{eqnarray}
u_{n} = a_n,
\
n\geq 0.
\label{pos}
\end{eqnarray}
\end{theorem}

\noindent
{\bf Proof:}
By Theorem 2.1, for $n\geq 0$,
\begin{eqnarray*}
G_{n}({\bf y})= {\cal K}^n  G_0({\bf y}).
\end{eqnarray*}
Putting ${\bf y} ={\bm \infty}$ gives (\ref{pos}).
$\Box$

For example,
\begin{eqnarray*}
u_0=a_0 = 1,
\
P(X_1\leq x)=u_1=a_1= E\ I(e_0\leq y_1) \int G_0(g_{\bf y} (sz_1+e_0), dz_1).
\end{eqnarray*}

\section{The case of $F$ absolutely continuous}
\setcounter{equation}{0}

Our solution Theorem 2.2 does not tell us how
$u_n$ behaves for large $n$.
Also calculating $a_n$ requires repeated integration.
Here, we give another solution that overcomes these problems,
using Fredholm integral theory given in Appendix A of Withers and Nadarajah (2009a), referred to below as the appendix.

\begin{theorem}
For ${\bf z} = (z_0,z_1)$ in $R^2$ and $h$ a function of ${\bf z}$, set $h_{.i}({\bf z})=\partial_i h({\bf z})$, where $\partial_i=\partial/\partial_i$.
Suppose that $F$ is absolutely continuous with pdf $f$ and that
$h({\bf z})\rightarrow 0$
as $z_1\rightarrow  \infty$.
Let $\delta(z_1)$ denote the Dirac delta function on $R$.
Set
\begin{eqnarray*}
&&
b_{\bf y} ({\bf z}) = f(y_{0x}-rz_0-sz_1),
\
\gamma_{\bf y} (z_0)=(y_{0x}-y_1-rz_0)/s,
\\
&&
A_{\bf yz} = \{rz_0+sz_1>y_{0x}-y_1\},
\\
&&
C_1({\bf y}, {\bf z}) = \delta(z_1-\gamma_{\bf y} (z_0)) b_{\bf y} ({\bf z}),
\
C_2({\bf y}, {\bf z}) = I(A_{\bf yz}) b_{{\bf y}.1} ({\bf z}),
\\
&&
C({\bf y}, {\bf z}) =\sum_{j=1}^2 C_j({\bf y}, {\bf z}),
\
K({\bf y}, {\bf z}) = -r C({\bf y}, {\bf z}).
\end{eqnarray*}
Then we can write (\ref{calK}) in the form
\begin{eqnarray}
{\cal K} r({\bf y}) = \int K({\bf y}, {\bf z}) r({\bf z}) d{\bf z}.
\label{kernel}
\end{eqnarray}
\end{theorem}

\noindent
{\bf Proof:}
Set $c_{\bf y}({\bf z}) = I(A_{\bf yz}) b_{\bf y}({\bf z})$. Then $c_{{\bf y}.1}({\bf z})=C({\bf y}, {\bf z})$.
Then
\begin{eqnarray*}
{\cal K}h({\bf y})=r\int\int dz_0 c_{\bf y}({\bf z}) h(z_0,dz_1)=-r\int h({\bf z}) c_{{\bf y}.1}({\bf z})d{\bf z},
\end{eqnarray*}
integrating by parts.
This ends the proof.
$\Box$

Although $K({\bf y}, {\bf z})$ is a generalized function, it satisfies

\begin{theorem}
For $r$, $|s|\neq 1$,
\begin{eqnarray*}
0< \int\int  K({\bf y}, {\bf z})K({\bf z}, {\bf y}) d{\bf y} d{\bf z} < \infty.
\end{eqnarray*}
\end{theorem}

\noindent
{\bf Proof:}
Note that
\begin{eqnarray*}
\int\int  C({\bf y}, {\bf z}) C({\bf z}, {\bf y}) d{\bf y} d{\bf z} = \sum_{i,j=1}^2 \alpha_{ij},
\end{eqnarray*}
where $\alpha_{ij}=\int\int  C_i({\bf y}, {\bf z}) C_j({\bf z}, {\bf y})d{\bf y} d {\bf z}$.
Note $\alpha_{11} $ involves two delta functions, so the four integrations
over $y_0$, $y_1$, $z_0$, $z_1$ reduce to two over  $y_0$, $z_0$ at $z_1=\gamma_{\bf y} (z_0)$, $y_1=\gamma_{\bf z} (y_0)$,
that is at $y_1=y_1^*$, $z_1=z_1^*$, where $y_1^*=[(y_{0x}-rz_0-s(z_{0x}-ry_0)]/(1-s^2)$
and $z_1^*=[(z_{0x}-ry_0-s(y_{0x}-rz_0)]/(1-s^2)$.
So,
\begin{eqnarray*}
\alpha_{11}=2\int\int_{y_0<z_0}  [K({\bf y}, {\bf z}) K({\bf z}, {\bf y})]_{y_1=y_1^*, z_1=z_1^*} dy_0dz_0=2(I_1+I_2+I_3),
\end{eqnarray*}
where $I_1$, $I_2$, $I_3$ integrate over $A_1=\{x<y_0<z_0\}$, $A_2=\{y_0<x<z_0\}$, and $A_3=\{y_0<z_0<x\}$.
A transformation of variables gives $I_i=b_i\int\int_{A_i} f({\bf u}) f({\bf v}) d{\bf u} d{\bf v}$,
where $b_1=1/b_2=(1-s^2)/r^2$ and $b_3=(1-s^{2})/(1-r^2)$.
The other  $\alpha_{ij}$ can be dealt with similarly.
$\Box$

This theorem implies that $K({\bf y}, {\bf z})$ is a (non-symmetric) Fredholm kernel with respect
to Lebesgue measure, allowing
the Fredholm theory of the appendix
to be applied, in particular the
functional forms of the Jordan form and singular value decomposition.

Let $\{\lambda_{j},r_{j},l_{j}:\ j\geq 1 \}$ be the eigenvalues and associated
right and left eigenfunctions of ${\cal K}$ ordered so that $|\lambda_{j}|\geq |\lambda_{j+1}|$.
If $\{\lambda_{j} \}$ are real then
$\{r_{j},l_{j} \}$ can be taken as real.
By the appendix referred to, these satisfy
\begin{eqnarray}
{\cal K} r_{j}({\bf y})=\lambda_{j} r_{j}({\bf y}),
\
\overline{l}_{j}({\bf z}) {\cal K}=\lambda_{j} \overline{l}_{j}({\bf z}),
\
\int r_{j}\overline{l}_{k} =\delta_{jk},
\label{lr}
\end{eqnarray}
where $\overline{\zeta}$ is the complex conjugate of $\zeta$,
$\overline{l}({\bf z}) {\cal K}=\int \overline{l}({\bf y}) K({\bf y}, {\bf z}) d{\bf y}$ and $\delta_{jk}$ is the Kronecker function.
So, $\{ r_{j}({\bf y}), l_{k}({\bf y}) \}$ are biorthogonal functions with respect to Lebesgue measure.

We now assume that $K({\bf y}, {\bf z})$ has diagonal Jordan form.
(This holds, for example, when the eigenvalues are distinct.
This will generally be the case for our applications.)
The functional equivalent of the Jordan form is, by (3.6) of Withers and Nadarajah (2008b),
\begin{eqnarray*}
K({\bf y}, {\bf z}) = \sum_{j=1}^\infty \lambda_j r_j({\bf y}) \overline{l}_j({\bf z}).
\end{eqnarray*}
This implies that
\begin{eqnarray}
K_n({\bf y}, {\bf z})={\cal K}^{n-1} K({\bf y}, {\bf z}) = \sum_{j=1}^\infty \lambda_j^n r_j({\bf y}) \overline{l}_j({\bf z}),
\label{it}
\end{eqnarray}
where ${\cal K}^n$ is the operator corresponding to the iterated kernel $K_n({\bf y}, {\bf z})$.
By (A.8) of Withers and Nadarajah (2009a) with $\mu$ Lebesgue measure on $R^2$, if ${\cal K} G$ is in $L_2(R^2)$ then
\begin{eqnarray}
{\cal K}^n G({\bf y}) = \sum_{j=1}^\infty B_j(G)\ r_j({\bf y}) \lambda_j^n,
\
n\geq 1,
\label{n}
\end{eqnarray}
where $B_j(G)=\int_{R^2} G\overline{l}_j$.
Putting ${\bf y} = \pmb{ \infty}$ and $G=G_0$ in (\ref{n}) gives

\begin{theorem}
For $B_j$ of (\ref{n}) and $n\geq 1$,
\begin{eqnarray*}
a_n=\sum_{j=1}^\infty r_j(\pmb{ \infty}) B_j(G_0) \lambda_j^n.
\end{eqnarray*}
\end{theorem}

\begin{corollary}
Suppose that the eigenvalue $\lambda_1$ of largest magnitude has multiplicity $M$.
For $n\geq 1$,
\begin{eqnarray*}
a_n  = B(G_0) \lambda_1^{n}(1+\epsilon_{n}),
\end{eqnarray*}
where $\epsilon_{n}\rightarrow 0$ exponentially as $n\rightarrow \infty$ and
\begin{eqnarray*}
B(G_0) =\sum_{j=1}^M  r_j(\pmb{\infty}) B_j(G_0).
\end{eqnarray*}
So, for $n\geq 1$, by (\ref{pos}), $u_{n} =  B(G_0) \lambda_1^{n}(1+\epsilon_{n})$.
\end{corollary}

Unfortunately, we cannot use the method of Withers and Nadarajah (2009c)
for the numerical solution of the equations for the eigenvalues and eigenfunctions.
For example, the first equation in (\ref{lr})
for $r({\bf y}) = r_j({\bf y})$ at $\lambda=\lambda_j$ can be written
$-\lambda r({\bf y})/r=c_1({\bf y})+c_2({\bf y})$, where $c_i({\bf y})=\int C_i ({\bf y}, {\bf z}) r({\bf z}) d{\bf z}$.
Suppose that we use Gaussian quadrature
\begin{eqnarray}
\int_{R^2} a({\bf z}) d{\bf z} \approx \sum_{j=1}^q w_j a({\bf z}_j),
\label{gauss}
\end{eqnarray}
where $\{{\bf z}_1, \cdots, {\bf z}_q \}$ are given points in  ${R^2}$ and $\{w_1,\cdots, w_q\}$ are given weights.
Then
\begin{eqnarray*}
c_2({\bf y}) \approx \sum_{j=1}^q w_j C_2({\bf y}, {\bf z}_j) r({\bf z}_j).
\end{eqnarray*}
However, $c_1({\bf y}) = \int r(z_0, \gamma_{\bf y} (z_0)) dz_0$ only has one
single integral, so would need a different approximation, say
\begin{eqnarray*}
\int_{R} a(z_0)dz_0 \approx \sum_{j=1}^{q'} w_j' a(z_j'),
\end{eqnarray*}
where $\{z_1',\cdots, z_{q'}'\}$ are given points in  ${R}$ and $\{w_1',\cdots, w_{q'}'\}$ are given weights.
This gives
\begin{eqnarray*}
c_1({\bf y}) \approx \sum_{j=1}^{q'} w_j' r(z_j', \gamma_{\bf y} (z_j')).
\end{eqnarray*}
Putting ${\bf y} = {\bf z}_j$ now gives the system of equations
\begin{eqnarray*}
-\lambda r({\bf z}_i)/r \approx \sum_{j=1}^{q'} w_j' r(z_j', \gamma_{{\bf z}_i}(z_j')) + \sum_{j=1}^q w_j C_2({\bf z}_i, {\bf z}_j) r({\bf z}_j),
\end{eqnarray*}
that is
\begin{eqnarray*}
-\lambda{\bf r}\approx {\bm \theta} + {\bf C} {\bf r},
\end{eqnarray*}
where ${\bf r}$ has $i$th element $r({\bf z}_i)/r$,
${\bf C}$ has $(i,j)$th element $w_j C_2({\bf z}_i, {\bf z}_j)$, but ${\bm \theta}$ has $i$th element
$\sum_{j=1}^{q'} w_j' r(z_j',\gamma_{{\bf z}_i}(z_j'))$, which is not a multiple of ${\bf r}$.

\begin{note}
Differentiating the first equation in (\ref{lr})
for $r({\bf y}) = r_j({\bf y})$ at $\lambda=\lambda_j$ gives the differential-integral
equation for $ r_{.1}({\bf y})$,
\begin{eqnarray*}
\lambda r_{.1}({\bf y})/f(y_1) = \int r_{.1}((y_{0x}-sz_1-y_1)/r, z_1)dz_1.
\end{eqnarray*}
\end{note}

\section{A numerical solution}
\setcounter{equation}{0}

In the last section we saw that the kernel  $K({\bf y}, {\bf z})$ is a generalized function and
that numerical solution by direct Gaussian quadrature fails.
Here, we show how to get around these problems by using the iterated Fredholm kernel $K_2({\bf y}, {\bf z})$ given by (\ref{it}).

By (\ref{kernel}),
\begin{eqnarray*}
K_2({\bf y}, {\bf z})=r^2\sum_{j,k=1}^2 C_{jk}({\bf y}, {\bf z}),
\end{eqnarray*}
where
\begin{eqnarray*}
C_{jk}({\bf y}, {\bf z})=\int C_j({\bf y}, {\bf t}) C_k ({\bf t}, {\bf z}) d{\bf t}.
\end{eqnarray*}
We first show that these are ordinary functions, not generalized functions.
This is clearly true for $j=k=2$.
Also $ C_{11}({\bf y}, {\bf z}) = b_{\bf y} ({\bf T}) b_{\bf T} ({\bf z})$
at ${\bf T} = {\bf T}({\bf y}, {\bf z})$ given by $T_1=\gamma_{\bf y} (T_0)$ and $z_1=\gamma_{\bf T} (z_0)$.
Eliminate $T_1$ using
\begin{eqnarray}
y_{0x}-y_1-rT_0=sT_1=s(T_{0x}-rz_0-sz_1).
\label{T1}
\end{eqnarray}
So, $rT_0+sT_{0x}=y_{0x}-y_1+s(rz_0+sz_1)=a$ say.
When $r>0$, $r+s>0$ this has a unique solution
\begin{eqnarray*}
T_0 = \min(a/(r+s),(a-sx)/r).
\end{eqnarray*}
Note that $T_1$ is then given by (\ref{T1}).
Also
\begin{eqnarray*}
&&
C_{12}({\bf y}, {\bf z}) = \int [b_{\bf y} ({\bf t}) C_2({\bf t}, {\bf z})]_{t_1=\gamma_{\bf y} (t_0)}dt_0,
\\
&&
C_{21}({\bf y}, {\bf z}) = \int [C_2 ({\bf y}, {\bf t}) b_{\bf t} ({\bf z})]_{z_1=\gamma_{\bf t} (z_0)}dt_0 =
\int [C_2({\bf y}, {\bf t}) b_{\bf t} ({\bf z})]_{t_1=t_{0x}-rz_0-sz_1}dt_0.
\end{eqnarray*}
This gives $K_2({\bf y}, {\bf z})$ as an ordinary function.
It follows that convergence in (\ref{it}) holds and Theorem 3.3 holds
for $n\geq 2$ although perhaps convergence does not hold at $n=1$.
Iterations can be done using:
\begin{eqnarray*}
a_{2n}= [{\cal K}^{2n}  G_0({\bf y})]_{{\bf y} = \pmb{\infty}},
\
a_{2n + 1}= [{\cal K}^{2n}  G_1({\bf y})]_{{\bf y} = \pmb{\infty}},
\
n\geq 0.
\end{eqnarray*}
Also $G_1({\bf y}) = {\cal K}G_0({\bf y})$ is an ordinary function since by (\ref{kernel})
\begin{eqnarray*}
{\cal K}G({\bf y}) =-r \sum_{j=1}^2 {\cal C}_j G({\bf y}),
\
{\cal C}_jG({\bf y})=\int C_j({\bf y}, {\bf z}) G({\bf z}) d{\bf z},
\end{eqnarray*}
and for $G({\bf y})$ an ordinary function, ${\cal C}_j G({\bf y})$ is also an ordinary function for $j=1,2$.
This is obvious for  $j=2$.
It is true for  $j=1$ since
\begin{eqnarray*}
{\cal C}_1G({\bf y}) = \int [b_{\bf y} ({\bf z}) G({\bf z})]_{z_1=\gamma_{\bf y} (z_0)} dz_0.
\end{eqnarray*}

Since
\begin{eqnarray*}
0<\int\int K_2({\bf y}, {\bf z}) K_2({\bf z}, {\bf y}) d{\bf y} d{\bf z} < \infty,
\end{eqnarray*}
(\ref{it}) holds for $n=2,4,\cdots$.
So, the left and right eigenfunctions  $r({\bf y}) = r_{j}({\bf y})$
and $\overline{l}({\bf y}) = \overline{l}_{j}({\bf y})$
for the eigenvalue $\theta=\lambda_{j}^2$ of $K_2({\bf y}, {\bf z})$ satisfy
\begin{eqnarray*}
{\cal K}^2 r({\bf y})=\theta r({\bf y}),
\
\overline{l}({\bf z}) {\cal K}^2=\theta \overline{l}({\bf z}),
\
\int r\overline{l}=1.
\end{eqnarray*}
By the Gaussian quadrature approximation, (\ref{gauss}), we can approximate this as
\begin{eqnarray*}
{\bf K}_2 {\bf r} \approx {\bm \theta} {\bf r},
\
{\bf K}_2' {\bf l} \approx \overline{\bm \theta} {\bf l},
\end{eqnarray*}
where ${\bf K}_2$ is the $q\times q$ matrix with $(i,j)$ element $w_jK_2({\bf z}_i, {\bf z}_j)$,
and ${\bf r}$ and ${\bf l}$ are the $q$-vectors with $j$th elements $w_j r({\bf z}_j)$ and $w_j l({\bf z}_j)$.
So, the first $q$ eigenvalues
and right and left eigenfunctions of $K_2({\bf y}, {\bf z})$ can be approximated by
the eigenvalues and right and left eigenvectors of ${\bf K}_2$ standardized so that they are biorthonormal.

So, finally we obtain the distribution of $M_n$ to be
\begin{eqnarray*}
u_{n} = a_{n}\approx \widehat{a}_{n},
\
n\geq 2,
\end{eqnarray*}
where
\begin{eqnarray*}
\widehat{a}_{2n+i} = \sum_{j=1}^q \widehat{\theta}_j^n  \widehat{r}_{jq} \widehat{B}_j(G_i),
\
n\geq 1,
\
i = 0, 1
\end{eqnarray*}
and
\begin{eqnarray*}
\widehat{B}_j(G_i) = \sum_{k=1}^q  w_k \overline{l}_{jk} G_i({\bf z}_k),
\end{eqnarray*}
where $\widehat{\theta}_j$ is the $j$th eigenvalue of  ${\bf K}_2$,
and $\widehat{r}_{jk}$ and $\overline{l}_{jk}$ are the $k$th components of its left and right eigenvectors.

For a more precise result one can let $q$ increase to $\infty$ with
$n$ and use known  expressions for the remainder in the Gaussian
approximation.
Compare equation (25.4.9) of Abramowitz and Stegun (1964).

\begin{note}
Using Theorem 2.2, Theorem 3.1, Theorem 3.3 and the numerical tools developed above
one can calculate $u_n = P (M_n < x)$ for absolutely continuous cdfs $F$ and $H$.
Figures 4.1 and 4.2 show plots of $u_n$ for $n = 1000$, $s = 1, 5$ and $r = 0.1, 0.2, \ldots, 0.8, 0.9$
when $G_0$ is a product of two independent standard normal cdfs.
In each figure, the distribution of $u_n$ becomes less dominant as $r$ increases from 0 to 1.
\end{note}

\newpage

\centerline{\epsfig{file=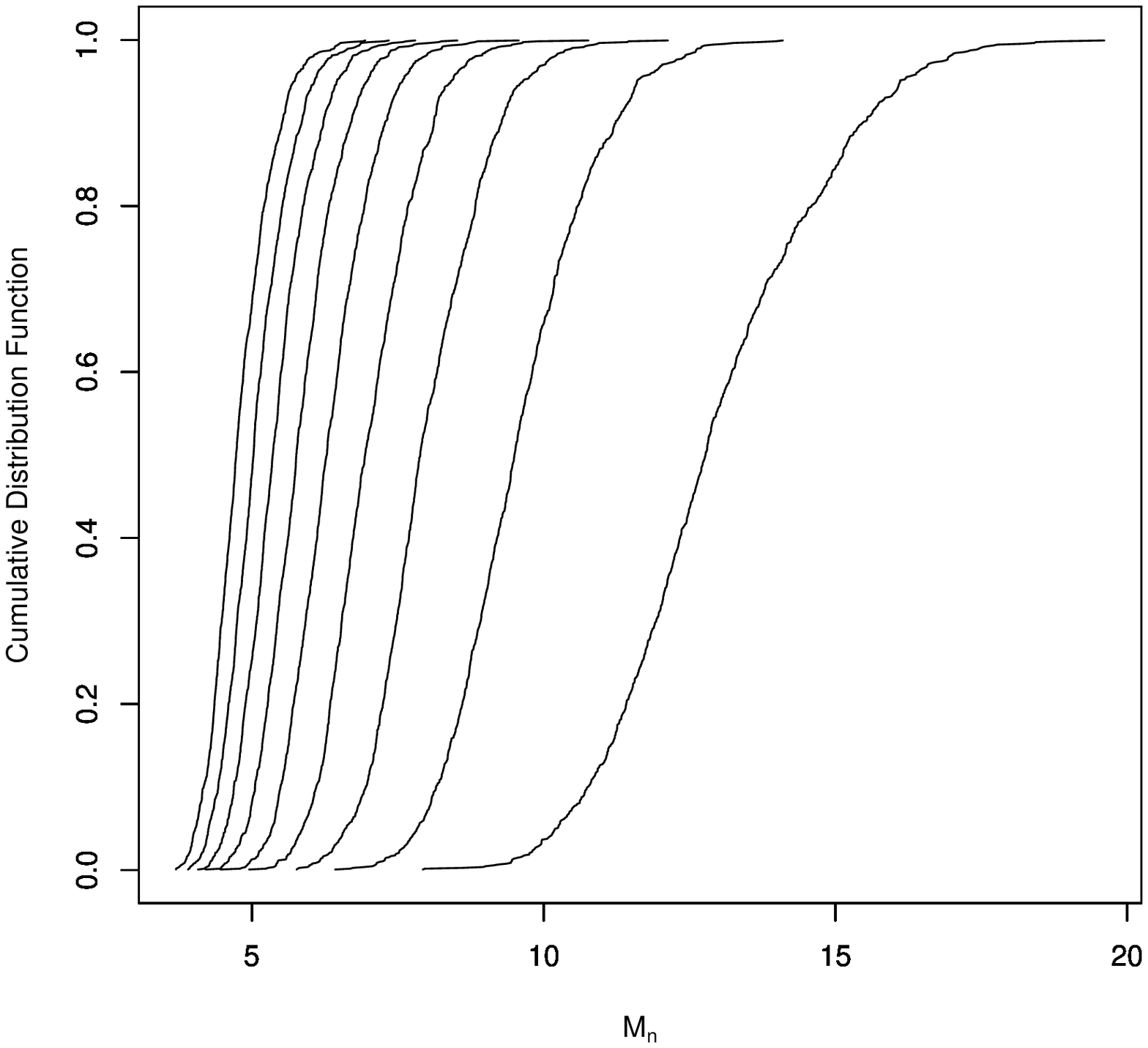,width=6in,height=3.8in}}
\noindent
{\bf Figure 4.1.}~~Plot of $u_n = P (M_n < x)$ versus $x$ for $n = 1000$, $s = 1$
and $r = 0.1$, $0.2$, $\ldots$, $0.8$, $0.9$
when $G_0$ is a product of two independent standard normal cdfs.
The curves from the left to right correspond to increasing values of $r$.

\centerline{\epsfig{file=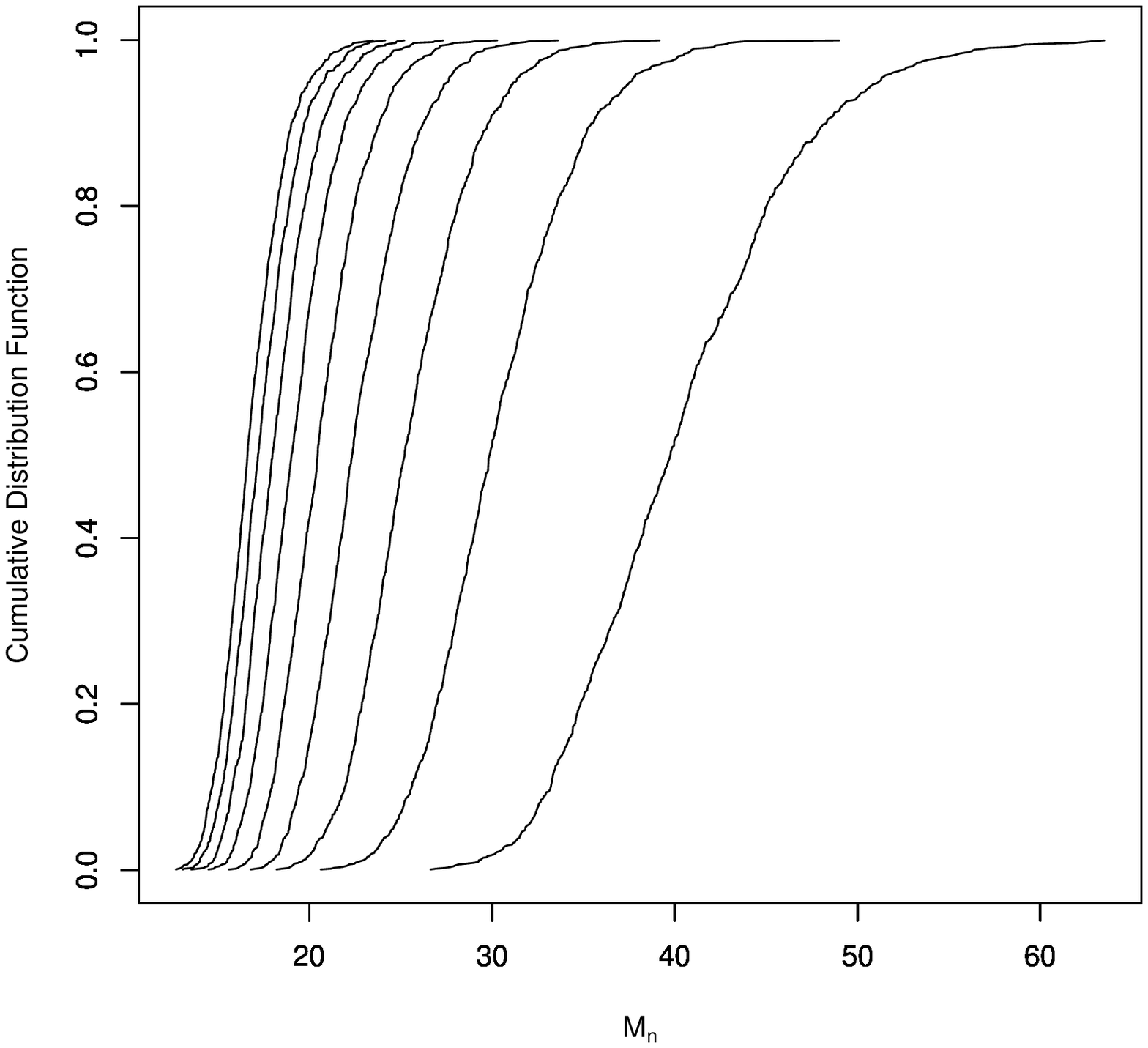,width=6in,height=3.8in}}
\noindent
{\bf Figure 4.2.}~~Plot of $u_n = P (M_n < x)$ versus $x$ for $n = 1000$, $s = 5$
and $r = 0.1$, $0.2$, $\ldots$, $0.8$, $0.9$
when $G_0$ is a product of two independent standard normal cdfs.
The curves from the left to right correspond to increasing values of $r$.

\end{document}